\documentclass[12pt,leqno]{article}
\usepackage{amsmath,amsthm,amscd,amssymb}
\usepackage{latexsym}
\newtheorem{remark}{Remark}[section]
\newcommand{\bremark}{\begin{remark} \em}
\newcommand{\eremark}{\end{remark} }

\begin{document}
\parindent 15pt
\renewcommand{\theequation}{\thesection.\arabic{equation}}
\renewcommand{\baselinestretch}{1.15}
\renewcommand{\arraystretch}{1.1}
\def\disp{\displaystyle}
\title{\bf\large A priori bounds for positive solutions of Kirchhoff type equations
\footnotetext{This work was supported by NNSF
of China(Grant: No.11271120)}
\footnotetext{E-mail: daiqiuyi@aliyun.com,\ Shifeilin1116@163.com}
\author{{\small Qiuyi Dai\ \ Enhao Lan\ \ Feilin Shi}\\
{\small Department of Mathematics, Hunan Normal University}\\
 {\small Changsha Hunan 410081, P.R.China}}
}
\maketitle

\abstract{Let $\Omega$ be a bounded smooth domain in $R^N$. Assume that $0<\alpha<\frac{2^*-1}{2}$, $a>0$, and $b>0$. We consider the following Dirichlet problem of Kirchhoff type equation
\begin{equation}
\left\{\begin{array}{ll}\label{1*}
-(a+b||\nabla u||_2^{2\alpha})\Delta u=|u|^{p-1}u+h(x,u,\nabla u), & \mbox{in}\ \ \Omega\\
u=0 & \mbox{on}\ \ \partial\Omega\\
\end{array}
\right.
\end{equation}
with $p\in (0,\ 2^*)\setminus\{1\}$. Where $2^*=+\infty$ for $N=2$, and $2^*=\frac{N+2}{N-2}$ for $N\geq 3$. Under suitable conditions of $h(x,u,\nabla u)$ (see $(A)$, $(H_1)$ and $(H_2)$ in section 3), we get a priori estimates for positive solutions to problem (\ref{1*}). By making use of these estimates and the continuous method, we further get existence results for positive solution to problem (\ref{1*}) when $0<p<1$, or $2\alpha+1<p<2^*$. Effects of the term $b||\nabla u||_2^{2\alpha+1}$ on the solution set of problem (\ref{1*}) can be seen from an example given in section 2.}

\section*{1 Introduction}

\setcounter{section}{1}

\setcounter{equation}{0}

\noindent Let $\Omega \subset R^N(N\geq 2)$ be an open bounded
smooth domain with boundary $\partial\Omega$, $||\bullet||_q$ denote the norm of $L^q(\Omega)$ for any $q>0$. Assume that $a>0$, $b>0$ and $0<\alpha<\frac{2^*-1}{2}$ with $2^*=+\infty$ for $N=2$, and $2^*=\frac{N+2}{N-2}$ for $N\geq 3$. We consider the following problem for Kirchhoff type equations
\begin{equation}
\left\{\begin{array}{ll}\label{1}
-(a+b||\nabla u||_2^{2\alpha})\Delta u=|u|^{p-1}u+h(x,u,\nabla u)& \mbox{in}\ \ \Omega,\\
u=0 & \mbox{on}\ \ \partial\Omega.
\end{array}
\right.
\end{equation}
When $\alpha=1$ and $h(x,u,\nabla u)=f(x,u)$, problem (\ref{1}) is reduced to the following problem involving standard Kirchhoff operator.
\begin{equation}
\left\{\begin{array}{ll}\label{12}
-(a+b||\nabla u||_2^{2})\Delta u=|u|^{p-1}u+f(x,u)& \mbox{in}\ \ \Omega,\\
u=0 & \mbox{on}\ \ \partial\Omega.
\end{array}
\right.
\end{equation}

In recent years, problems like (\ref{12}) have been extensively studied by making use of variational method. Mention some but few, we refer to \cite{sAC,sCK,sCL,sCW,sLLS1,sLLT,sM,sMR,sMZ,sPZ,sST,sYZ,sZP, sA,sHZ,sJW,sLG,sLLS,sNW, sWTXZ,sX, sDN, sFig, sFS} and the references cited there in.

In general, problem (\ref{1}) has no variational structure. Hence, we adopt a device of using a priori estimate and fixed point theorem to study the positive solutions to it. This device has its benefits for it can ensure compactness of  the positive solution set of problem (\ref{1}) \cite{sFLN, sGS, sFY, sRu}. To make this device work, a crucial step is to drive an a priori estimate for solutions to the problem under consideration. It is believed that the existence result can be ensured once a priori estimates of solutions are established. However, this confidence is violated in our model problem (\ref{1}). In fact, for arbitrary $a, b>0$, we can get easily upper and lower bounds for nontrivial solutions to problem (\ref{1}) when $1<p<2\alpha+1$ and $h(x,u,\nabla u)\equiv 0$. Whereas, we find that problem (\ref{1}) has nontrivial solution only for some $a$ and $b$ in this case (see Theorem 2.1(ii) in section 2).

The rest of this paper is arranged as the following. Section 2 is used to analysis the solvability of problem (\ref{1}) when $h(x,u,\nabla u)\equiv 0$. Section 3 devotes to derive a priori estimates of positive solutions to (\ref{1}). Existence results for positive solutions to problem (\ref{1}) with $0<p<1$, or $2\alpha+1<p<2^*$ are given in section 4.

\section*{2 A simple observation for $h(x,u,\nabla u)\equiv 0$}

\setcounter{section}{2}

\setcounter{equation}{0}
Assume that $p\in(0,\ 2^*)\setminus\{1\}$ and $0<\alpha<\frac{2^*-1}{2}$. In this section, we analysis the solvability of the following problem
\begin{equation}\label{eq151*}
\left\{
\begin{array}{ll}
-(a+b||\nabla u||_2^{2\alpha})\Delta u=|u|^{p-1}u,&x\in\Omega,\\
u=0,&x\in\partial\Omega.
\end{array}
\right.
\end{equation}

It is easy to see that if $u$ is a nontrivial solution of problem (\ref{eq151*}) and $v=\eta u$ with $\eta=(a+b||\nabla u||_2^{2\alpha})^{\frac{1}{1-p}}$, then $v$ is a nontrivial solution of the following well studied problem
\begin{equation}\label{151*}
\left\{
\begin{array}{ll}
-\Delta v=|v|^{p-1}v,& x\in\Omega\\
v=0,& x\in\partial\Omega
\end{array}
\right.
\end{equation}

Hence, any nontrivial solution of problem (\ref{eq151*}) can be obtained as the form $u=(a+b\beta^{\alpha})^{\frac{1}{p-1}} v$ with $v$ being a nontrivial solution of problem (\ref{151*}) and $\beta$ being a positive solution of the following algebraic equation
$$y^{\frac{p-1}{2}}-b||\nabla v||_2^{p-1}y^\alpha-a||\nabla v||_2^{p-1}=0.$$

Let
$$S(\Omega)=\inf\limits_{0\not\equiv u\in H^1_0(\Omega)}\frac{||\nabla u||_2^2}{||u||^2_{p+1}}.$$

Since $0<p<2^*$, it is well known that $S(\Omega)$ can be attained by some positive function $v\in H^1_0(\Omega)$. Let $v$ be a solution of  problem (\ref{151*}) which attained $S(\Omega)$. Then, it is easy to check that
$$||\nabla v||_2^{p-1}=S^{\frac{p+1}{2}}(\Omega).$$
For the simplicity of the notations, we set $S=S^{\frac{p+1}{2}}(\Omega)$ and $\gamma=2\alpha+1-p$. Let
$$f(y)=y^{\frac{p-1}{2}}-bSy^\alpha-aS.$$
For $y>0$, it is easy to see that
\begin{equation}\label{152*}
\left\{
\begin{array}{ll}
f^\prime(y)=y^{\alpha-1}(\frac{p-1}{2}y^{-\frac{\gamma}{2}}-\alpha bS)\\
f^{\prime\prime}(y)=y^{\alpha-2}(\frac{(p-1)(p-3)}{4}y^{-\frac{\gamma}{2}}-\alpha(\alpha-1)bS)
\end{array}
\right.
\end{equation}
Therefore, we have

(i)\ If $0<p<1$, then $f^\prime(y)<0$ for $y\in(0,\ +\infty)$, $\lim\limits_{y\rightarrow 0}f(y)=+\infty$, and $\lim\limits_{y\rightarrow +\infty}f(y)=-\infty$. Hence, the equation $f(y)=0$ has a unique solution in $(0,\ +\infty)$.

(ii)\ If $1<p<2\alpha+1$, then $\gamma>0$, and the equation $f^\prime(y)=0$ has a unique solution $y_0=(\frac{p-1}{2\alpha bS})^{\frac{2}{\gamma}}$ in $(0,\ +\infty)$. Moreover,
\begin{equation}\label{153*}
\left\{
\begin{array}{ll}
f^\prime(y)>0, &\mbox{for} \ \ y<y_0\\
f^\prime(y)<0, &\mbox{for} \ \ y>y_0\\
\max\limits_{y\in(0,\ +\infty)}\{f(y)\}=\gamma[\frac{(p-1)^{(p-1)}}{(2\alpha)^{(2\alpha)}}]^{\frac{1}{\gamma}}(bS)^{\frac{1-p}{\gamma}}-aS
\end{array}
\right.
\end{equation}
From this, we get
\begin{equation}\label{154*}
\left\{
\begin{array}{ll}
\max\limits_{y\in(0,\ +\infty)}\{f(y)\}>0,& \mbox{for}\ \ ab^{\frac{p-1}{\gamma}}<\gamma[\frac{(p-1)^{(p-1)}}{(2\alpha S)^{(2\alpha)}}]^{\frac{1}{\gamma}}\\
\max\limits_{y\in(0,\ +\infty)}\{f(y)\}=0,& \mbox{for}\ \ ab^{\frac{p-1}{\gamma}}=\gamma[\frac{(p-1)^{(p-1)}}{(2\alpha S)^{(2\alpha)}}]^{\frac{1}{\gamma}}\\
\max\limits_{y\in(0,\ +\infty)}\{f(y)\}<0,& \mbox{for}\ \ ab^{\frac{p-1}{\gamma}}>\gamma[\frac{(p-1)^{(p-1)}}{(2\alpha S)^{(2\alpha)}}]^{\frac{1}{\gamma}}\\
\end{array}
\right.
\end{equation}
Noting that $\lim\limits_{y\rightarrow 0}f(y)=-aS<0$, and $\lim\limits_{y\rightarrow +\infty}f(y)=-\infty$, we can conclude from (\ref{154*}) that the equation $f(y)=0$ has
\begin{equation}\label{155*}
\left\{
\begin{array}{ll}
\mbox{exactly two solutions in}\ (0,\ +\infty), & \mbox{when}\ \ ab^{\frac{p-1}{\gamma}}<\gamma[\frac{(p-1)^{(p-1)}}{(2\alpha S)^{(2\alpha)}}]^{\frac{1}{\gamma}},\\
\mbox{exactly one solution in}\ (0,\ +\infty), & \mbox{when}\ \ ab^{\frac{p-1}{\gamma}}=\gamma[\frac{(p-1)^{(p-1)}}{(2\alpha S)^{(2\alpha)}}]^{\frac{1}{\gamma}},\\
\mbox{no solution in}\ (0,\ +\infty), & \mbox{when}\ \ ab^{\frac{p-1}{\gamma}}>\gamma[\frac{(p-1)^{(p-1)}}{(2\alpha S)^{(2\alpha)}}]^{\frac{1}{\gamma}}.
\end{array}
\right.
\end{equation}

In this case, it is easy to see that the critical point $y_0$ of $f(y)$ satisfies
\begin{equation}\label{ad1}
\lim\limits_{b\rightarrow 0}y_0=+\infty.
\end{equation}
Moreover, if $ab^{\frac{p-1}{\gamma}}<\gamma[\frac{(p-1)^{(p-1)}}{(2\alpha S)^{(2\alpha)}}]^{\frac{1}{\gamma}}$ and $0<y_1<y_2$ are two solutions of the equation $f(y)=0$, then we have
\begin{equation}\label{ad2}
\lim\limits_{b\rightarrow 0}y_1=(aS)^{\frac{2}{p-1}},\;\;\;\lim\limits_{b\rightarrow 0}y_2=+\infty.
\end{equation}

(iii)\ If $p>2\alpha+1$, then $\gamma<0$, and the equation $f^\prime(y)=0$ has a unique solution $y_0=(\frac{p-1}{2\alpha bS})^{\frac{2}{\gamma}}$ in $(0,\ +\infty)$. Moreover,
\begin{equation}\label{156*}
\left\{
\begin{array}{ll}
f^\prime(y)<0, &\mbox{for} \ \ y<y_0\\
f^\prime(y)>0, &\mbox{for} \ \ y>y_0\\
\min\limits_{y\in(0,\ +\infty)}\{f(y)\}=\gamma[\frac{(p-1)^{(p-1)}}{(2\alpha)^{(2\alpha)}}]^{\frac{1}{\gamma}}(bS)^{\frac{1-p}{\gamma}}-aS
\end{array}
\right.
\end{equation}
Noting that $\lim\limits_{y\rightarrow 0}f(y)=-aS<0$, and $\lim\limits_{y\rightarrow +\infty}f(y)=+\infty$, we can conclude from (\ref{156*}) that the equation $f(y)=0$ has exactly one solution in $(0,\ +\infty)$.

(iv)\ If $p=2\alpha+1$, then it is easy to see that the equation $f(y)=0$ has exactly one solution in $(0,\ +\infty)$ when $bS<1$, and has no solution in $(0,\ +\infty)$ when $bS\geq 1$.

From the above discussions, we finally reach the following conclusions.
\vskip 0.2cm
{\em {\bf Theorem 2.1.}\ Assume that $a,\ b>0$, $\alpha>0$, and $p\in(0,\ 2^*)\setminus\{1\}$. then the following claims hold

(i)\ If $0<p<1$, or $p>2\alpha+1$, then problem (\ref{eq151*}) has as many solutions as problem (\ref{151*}).

(ii)\ If $1<p<2\alpha+1$, then problem (\ref{eq151*}) has at least two positive solutions when $ab^{\frac{p-1}{\gamma}}<\gamma[\frac{(p-1)^{(p-1)}}{(2\alpha S)^{(2\alpha)}}]^{\frac{1}{\gamma}}$, has at least one positive solution when $ab^{\frac{p-1}{\gamma}}=\gamma[\frac{(p-1)^{(p-1)}}{(2\alpha S)^{(2\alpha)}}]^{\frac{1}{\gamma}}$, and has no nontrivial solution when $ab^{\frac{p-1}{\gamma}}>\gamma[\frac{(p-1)^{(p-1)}}{(2\alpha S)^{(2\alpha)}}]^{\frac{1}{\gamma}}$.

(iii)\ If $p=2\alpha+1$, then problem (\ref{eq151*}) has at least one positive solution when $bS<1$, and has no nontrivial solution when $bS\geq 1$.}
\vskip 0.1cm
From (\ref{ad1}) and (\ref{ad2}), we can see that
\vskip 0.1cm
{\em {\bf Theorem 2.2.}\ Assume that $a,\ b>0$, $\alpha>0$, and $1<p<2\alpha+1$. If $ab^{\frac{p-1}{\gamma}}\leq\gamma[\frac{(p-1)^{(p-1)}}{(2\alpha S)^{(2\alpha)}}]^{\frac{1}{\gamma}}$, then problem (\ref{eq151*}) has at least one positive solution $u_b(x)$ such that
$$\lim\limits_{b\rightarrow 0}\max\limits_{x\in\Omega}u_b(x)=+\infty.$$}

Contrast to the case $1<p<2\alpha+1$, we have
\vskip 0.1cm
{\em {\bf Theorem 2.3.}\ Assume that $a,\ b>0$, $\alpha>0$. If $0<p<1$, or $2\alpha+1<p<2^*$, then any positive solution $u_b(x)$ of problem (\ref{eq151*}) converges to a positive solution of the following problem
\begin{equation}\label{ad3}
\left\{
\begin{array}{ll}
-a\Delta v=|v|^{p-1}v,& x\in\Omega\\
v=0,& x\in\partial\Omega,
\end{array}
\right.
\end{equation}
when $b\rightarrow 0$.}

The conclusion of Theorem 2.3 can be seen from the uniform a priori estimates given in Theorem 3.7 and 3.8 of section 3.

\section*{3 A priori estimates}

\setcounter{section}{3}

\setcounter{equation}{0}

This section devotes to derive a uniform bound for the solutions to
problem (\ref{1}). For the need of proving existence results, we consider the following problem of Kirchhoff type equations with a parameter $t\in[0,1]$.

\begin{equation}
\left\{\begin{array}{ll}\label{eq1}
-A(t,u)\Delta u=u^p+th(x,u,\nabla u), & \mbox{in}\ \ \Omega,\\
u>0 & \mbox{in}\ \ \Omega,\\
u=0 & \mbox{on}\ \ \partial\Omega,\\
\end{array}
\right.
\end{equation}
where $a>0$, $b>0$, $A(t,u)=a+tb||\nabla u||_2^{2\alpha}$, and $h(x,s,\xi)$ satisfies
\vskip 0.1cm
{\bf $(A)$.}\ For $s>0$, $h(x,s,\xi)$ is Holder continuous and $h(x,s,\xi)\geq 0$ .
\vskip 0.1cm
{\bf $(H_1)$.}\ If $0<p<1$, we assume that $\lim\limits_{s\rightarrow 0}\frac{h(x,s,\xi)}{s^p}=0$ uniformly in $x$ and $\xi$, and $|h(x,s,\xi)|\leq C(s^q+s^{q_1})$ for $s>0$ and $0<p<q\leq q_1<1$.
\vskip 0.1cm
{\bf $(H_2)$.}\ If $1<p<2^*$, we assume that there exists a positive constant $\lambda$ such that $\lim\limits_{s\rightarrow 0}\frac{h(x,s,\xi)}{s}=\lambda$ uniformly in $x$ and $\xi$, and $|h(x,s,\xi)|\leq C(s+s^q)$ for $s>0$ and $1<q<p<2^*$.

\vskip 0.1cm
{\em {\bf Remark 3.1.}\ If $0<p<1$, then $h(x,u,\nabla u)=\mu |u|^{q-1}u+\frac{|u|^{q_1-1}u|\nabla u|^2}{1+|\nabla u|^2}$ satisfies $(A)$ and $(H_1)$. If $1<p<2^*$, then $h(x,u,\nabla u)=\lambda u+\frac{|u|^{q-1}u|\nabla u|^2}{1+|\nabla u|^2}$ satisfies $(A)$ and $(H_2)$.}
\vskip 0.1cm
{\em {\bf Remark 3.2.}\ All our results in this section maybe hold under assumptions weaker than $(H_1)$ and $(H_2)$ (see \cite{sFY, sRu}). However, for the simplicity, we only prove our results under the assumptions $(H_1)$ and $(H_2)$.}
\vskip 0.1cm
To get upper bound of positive solutions to (\ref{eq1}) for $1<p<2^*$, we need following two well known results for semilinear problem.
\vskip 0.2cm
{\em {\bf Lemma 3.3.(\cite{sGS})} If $1<p<2^*$, then the following problem has no nontrivial $C^2$ solution.
\begin{equation}
\left\{\begin{array}{ll}\label{13}
-\Delta u=|u|^{p-1}u & \mbox{in}\ \ R^N,\\
u\geq 0 & \mbox{in}\ \ R^N.
\end{array}
\right.
\end{equation}}
\vskip 0.2cm
{\em {\bf Lemma 3.4.(\cite{sGS})} If $1<p<2^*$, then the following problem has no nontrivial $C^2$ solution.
\begin{equation}
\left\{\begin{array}{ll}\label{14}
-\Delta u=|u|^{p-1}u & \mbox{in}\ \ R^N_c=\{(x_1,\cdots,x_N)\in R^N\|\ \ x_N>c\ \},\\
u\geq 0 & \mbox{in}\ \ R^N_c,\\
u=0 & x_N=c.
\end{array}
\right.
\end{equation}}
\vskip 0.2cm
To get lower bound of positive solutions to (\ref{eq1}) for $0<p<1$, we need following results.
\vskip 0.2cm
{\em {\bf Lemma 3.5.} If $0<p<1$, then the following problem has no bounded nontrivial $C^2$ solution.
\begin{equation}
\left\{\begin{array}{ll}\label{15}
-\Delta u=|u|^{p-1}u & \mbox{in}\ \ R^N,\\
u\geq 0 & \mbox{in}\ \ R^N.
\end{array}
\right.
\end{equation}}
\vskip 0.2cm
{\em {\bf Lemma 3.6.} If $0<p<1$, then the following problem has no bounded nontrivial $C^2$ solution.
\begin{equation}
\left\{\begin{array}{ll}\label{16}
-\Delta u=|u|^{p-1}u & \mbox{in}\ \ R^N_c,\\
u\geq 0 & \mbox{in}\ \ R^N_c,\\
u=0 & x_N=c.
\end{array}
\right.
\end{equation}}
\vskip 0.2cm
Lemma 3.5 and 3.6 maybe also well known, but we can not fined a suitable reference, so we give a proof of Lemma 3.5 here. The proof of Lemma 3.6 is very similar.
\vskip 0.2cm
{\bf Proof of Lemma 3.5.}\ Suppose by contradiction that problem (\ref{15}) has a bounded nontrivial solution $u$. Then by strong maximum principle we have $u(x)>0$ for any $x\in R^N$. Let $M=\sup\limits_{x\in R^N}u(x)$, $v(x)=u(x)/M$ and $h(y)=v(M^{\frac{1-p}{2}}y)$. Then $h(y)$ satisfies
\begin{equation}
\left\{\begin{array}{ll}\label{17}
-\Delta h=h^p & y\in R^N,\\
0<h\leq 1 & y\in R^N.
\end{array}
\right.
\end{equation}

Let $B_R(y_0)$ denote the ball in $R^N$ with center $y_0$ and radius $R$. Let $\lambda_1(R)$ denote the first eigenvalue of the following eigenvalue problem
\begin{equation}
\left\{\begin{array}{ll}\label{18}
-\Delta\varphi=\lambda\varphi & y\in B_R(y_0),\\
\varphi=0 & y\in\partial B_R(y_0)
\end{array}
\right.
\end{equation}
We choose $R$ so large that $\lambda_1(R)<1$. Denoting by $\varphi_1(y)>0$ the eigenfunction corresponding to $\lambda_1(R)$, and by $n$ the unit outward norm vector field on $\partial B_R(y_0)$, then we have $\frac{\partial\varphi_1}{\partial n}<0$ on $\partial B_R(y_0)$. Multiplying the differential equation in (\ref{17}) by $\varphi_1(y)$ and integrating the result equation on $B_R(y_0)$, we get
$$\int_{\partial B_R(y_0)}h\frac{\partial\varphi_1}{\partial n}dS=\int_{B_R(y_0)}(h^{p-1}-\lambda_1(R))h\varphi_1dy.$$
Since $0<h\leq 1$ on $\overline{B_R(y_0)}$, and $0<p<1$, we have $h^{p-1}-\lambda_1(R)>0$ on $B_R(y_0)$.

Therefore, on one hand, we have
$$\int_{B_R(y_0)}(h^{p-1}-\lambda_1(R))h\varphi_1dy>0.$$
On the other hand, we have
$$\int_{\partial B_R(y_0)}h\frac{\partial\varphi_1}{\partial n}dS<0,$$
due to $h>0$ and $\frac{\partial\varphi_1}{\partial n}<0$ on $\partial B_R(y_0)$. This reaches a contradiction.

At this stage, we are ready to prove our main results of this section
\vskip 0.1cm
{\em {\bf Theorem 3.7.} Assume that $0<p<1$, and that $(A)$ and $(H_1)$ hold. Then there exist positive constants $r$ and $R$ which are independent of $t$ and the solution of problem (\ref{eq1}) such that for any solution \mbox{$u\in C^2(\Omega)\cap C(\overline{\Omega})$} of problem (\ref{eq1}), there holds $r\leq||u||_{C(\Omega)}\leq R$.}
\vskip 0.1cm
{\bf Proof.}\  Multiplying the differential equation in problem (\ref{eq1}) by $u$ and integrating the result equation on $\Omega$, we obtain
\begin{equation}\label{eq151}
A(t,u)||\nabla u||_2^2=||u||_{p+1}^{p+1}+t\int_{\Omega}uh(x,u,\nabla u)dx
\end{equation}
Taking the assumption $(H_1)$ into account, we can deduce from (\ref{eq151}) that
\begin{equation}\label{eq19}
A(t,u)||\nabla u||_2^2\leq||u||_{p+1}^{p+1}+C(||u||_{q+1}^{q+1}+||u||_{q_1+1}^{q_1+1}).
\end{equation}
Since $0<p<q\leq q_1<1$ and
\begin{equation}\label{eq152}
A(t,u)||\nabla u||_2^2\geq a||\nabla u||_2^2 \ \ \ \mbox{for}\ t\in [0,\ 1],
\end{equation}
it follows from (\ref{eq19}) and the Holder inequality that
\begin{equation}\label{eq153}
||\nabla u||_2^2\leq C_1||u||_{q_1+1}^{q_1+1}+C_2 \ \ \  \mbox{for}\ t\in [0,\ 1]
\end{equation}

By the Sobolev embedding theorem, we know that there exists a positive constant $S_1(\Omega)$ which is independent of $u$ such that
\begin{equation}\label{eq154}
S_1(\Omega)||\nabla u||_2^{q_1+1}\geq ||u||^{q_1+1}_{q_1+1}.
\end{equation}
Since $0<q_1<1$, we can deduce from (\ref{eq153}) and (\ref{eq154}) that $||\nabla u||_{2}\leq C$ with $C$ being a universal positive constant. Now, a standard bootstrap argument will imply that there is a universal positive constant $R$ such that $||u||_{C(\Omega)}\leq R$.

To get a lower bound, we adopt a contradiction argument. Suppose by contradiction that there are sequences $\{t_n\}_{n=1}^{\infty}\subset(0,\ 1]$, $\{u_n\}_{n=1}^{\infty}$ and $\{x_n\}_{n=1}^{\infty}\subset\Omega$ such that
\begin{equation}
\left\{\begin{array}{ll}\label{eq20}
-A(t_n,u_n)\Delta u_n=u_n^p+t_nh(x,u_n,\nabla u_n), & \mbox{in}\ \ \Omega,\\
u_n>0 & \mbox{in}\ \ \Omega,\\
u_n=0 & \mbox{on}\ \ \partial\Omega,\\
\end{array}
\right.
\end{equation}
and $M_n=u_n(x_n)=\max\limits_{x\in\Omega}u_n(x)\rightarrow 0$ as $n\rightarrow\infty$. Since $\{x_n\}$ is bounded, up to a subsequence, we may assume that $x_n\rightarrow x_0\in\overline{\Omega}$ as $n\rightarrow\infty$. By (\ref{eq19}) and (\ref{eq152}), we have
$$||\nabla u_n||_2^2\leq C(M_n^{p+1}+M_n^{q+1}+M_n^{q_1+1})\rightarrow 0\ \ \mbox{as}\ \ n\rightarrow\infty.$$
Hence, $A(t_n,u_n)\rightarrow a>0$ as $n\rightarrow+\infty$. Set $\rho_n=M_n^{\frac{p-1}{2}}/\sqrt{A(t_n,u_n)}$. Then we have
$$\lim\limits_{n\rightarrow+\infty}\rho_n=+\infty,$$
due to $0<p<1$ and $M_n\rightarrow 0$ as $n\rightarrow+\infty$.

Let
$$y=\rho_n(x-x_n),\ \ \ \ \Omega_n=\rho_n(\Omega-\{x_n\}),\ \ \ \ v_n(y)=\frac{u(x_n+\rho_n^{-1}y)}{M_n}.$$
Then, $v_n(y)$ satisfies
\begin{equation}
\left\{\begin{array}{ll}\label{eq21}
-\Delta v_n=v_n^p+t_n\frac{h(x,M_nv_n,M_n\nabla v_n)}{M_n^p}, & y\in\Omega_n,\\
v_n>0 & y\in\Omega_n,\\
v_n=0 & y\in\partial\Omega_n,\\
\end{array}
\right.
\end{equation}

Now, we have the following two cases to be considered.

{\bf Case I.}\ If $\{\mbox{dist}(x_n,\partial\Omega_n)\}_{n=1}^\infty$ is unbounded, then, up to a subsequence, we may assume that $\mbox{dist}(x_n,\partial\Omega_n)\rightarrow\infty$ and $\Omega_n\rightarrow R^N$ as $n\rightarrow\infty$. Therefore, for any fixed compact domain $K\subset R^N$ there exist $n_K$ such that $K\subset\Omega_n$ when $n>n_K$. Since $0<v_n\leq 1$, we can get from the assumption $(H_1)$ that
$$|v_n^p+t_n\frac{h(x,M_nv_n,M_n\nabla v_n)}{M_n^p}|\leq C,$$
for some universal positive constant $C$. By the standard regularity theory of elliptic equations, we know that there exists a positive constant $C$ independent of $n$ such that for some $\tau\in(0,1)$ there holds
$$||v_n||_{C^{1,\tau}(K)}\leq C,\ \ \ \mbox{for}\ \ \ \ n>n_K.$$
Hence, up to a subsequence, we have $v_n$ converges in $C^1(K)$ to some function $v$. Since $K$ is arbitrary, by a diagonal process, we can choose a subsequence $\{v_{n_k}\}$ of $\{v_n\}$ such that $v_{n_k}\rightarrow v$ in $C^1_{loc}(R^N)$ as $k\rightarrow +\infty$. Taking the assumption $(H_1)$ into account, it is easy to see that $v$ satisfies
\begin{equation}
\left\{\begin{array}{ll}\label{22}
-\Delta v=v^p & \mbox{in}\ \ R^N,\\
v(0)=1,\\
v\geq 0 & \mbox{in}\ \ R^N.
\end{array}
\right.
\end{equation}
This contradicts with the conclusion of Lemma 3.5.

{\bf Case II.}\ If $\{\mbox{dist}(x_n,\partial\Omega_n)\}_{n=1}^\infty$ is bounded, then, up to a subsequence, we may assume that $\mbox{dist}(x_n,\partial\Omega_n)\rightarrow c$ for some constant $c$, and $\Omega_n\rightarrow R^N_{-c}$ as $n\rightarrow\infty$. Similar to the Case I, we can find a nonnegative function $v$ such that
\begin{equation}
\left\{\begin{array}{ll}\label{23}
-\Delta v=v^p & \mbox{in}\ \ R^N_{-c},\\
v(0)=1,\\
u=0 & x_N=-c.
\end{array}
\right.
\end{equation}
This contradicts with the conclusion of Lemma 3.6.

Therefore, there exists a universal positive constant $r$ such that for any solution $u$ of (\ref{eq1}), there hold $||u||_{C(\Omega)}\geq r$.

\vskip 0.1cm
Let $\lambda_1(\Omega)$ denote the first eigenvalue of the following eigenvalue problem
\begin{equation}
\left\{\begin{array}{ll}\label{24}
-\Delta\varphi=\lambda\varphi & y\in \Omega,\\
\varphi=0 & y\in\partial \Omega.
\end{array}
\right.
\end{equation}
Then, we have
\vskip 0.1cm
{\em {\bf Theorem 3.8.}\ Assume that $2\alpha+1<p<2^*$, and that $(A)$ and $(H_2)$ hold with $\lambda<a\lambda_1(\Omega)$. Then there exist two universal positive constants $r$ and $R$ such that for any solution \mbox{$u\in C^2(\Omega)\cap C(\overline{\Omega})$} of problem (\ref{eq1}), there holds $r\leq||u||_{C(\Omega)}\leq R$.}

{\bf Proof.}\ We prove the Theorem by contradiction. To get the lower bound, we suppose that there exist a sequence $\{t_n\}\subset(0,\ 1]$ and a corresponding solution sequence $\{u_n\}$ of problem (\ref{eq1}) such that
$$M_n=||u_n||_{C(\Omega)}\rightarrow 0,\ \ \mbox{and}\ \ t_n\rightarrow t_0,\ \ \mbox{as}\ \ n\rightarrow +\infty.$$
By the differential equation in (\ref{eq1}) and the assumption $(H_2)$, we can deduce that
$$||\nabla u_n||^2_2\leq ||u||_{p+1}^{p+1}+C(||u||_2^2+||u||_{q+1}^{q+1})\leq CM_n^2\rightarrow 0,\ \ \mbox{as}\ \ n\rightarrow+\infty.$$

Let $u_n(x)=M_nv_n(x)$. Then $v_n(x)$ solve the following problem
\begin{equation}\label{eq156}
\left\{
\begin{array}{ll}
-A(t_n,u_n)\Delta v_n=M_n^{p-1}v_n^p+t_n\frac{h(x,M_nv_n,M_n\nabla v_n)}{M_n},&x\in\Omega,\\
v_n>0,&x\in\Omega,\\
v_n=0,&x\in\partial\Omega.
\end{array}
\right.
\end{equation}
By $(H_2)$ and the standard elliptic estimates, we can easily see that, up to a subsequence, $v_n$ converges in $C^2(\Omega)$ to a positive function $v$. Moreover $v$ satisfies
\begin{equation}\label{eq1567}
\left\{
\begin{array}{ll}
-a\Delta v=t_0\lambda v,&x\in\Omega,\\
v=0,&x\in\partial\Omega.
\end{array}
\right.
\end{equation}
On the other hand, problem (\ref{eq1567}) has no positive solution due to $t_0\lambda\leq\lambda<a\lambda_1(\Omega)$. This reaches a contradiction. Therefore, there is a universal constants $r>0$ such that for any positive solution of (\ref{eq1}) there holds $||u||_{C(\Omega)}\geq r$.

To get the upper bound, we suppose that there exist a sequence $\{t_n\}\subset(0,\ 1]$, a corresponding solution sequence $\{u_n\}$ of problem (\ref{eq1}), and a sequence $\{x_n\}$ such that
$$M_n=u_n(x_n)=||u_n||_{C(\Omega)}\rightarrow+\infty,\ \ \ t_n\rightarrow t_0,\ \ \mbox{and}\ \ \ x_n\rightarrow x_0\in\overline{\Omega},\ \ \mbox{as}\ \ n\rightarrow +\infty.$$

Let $u_n(x)=M_nv_n(x)$, $y=\rho_n(x-x_n)$, $\Omega_n=\rho_n(\Omega-\{x_n\})$. If we choose $\rho_n$ so that
$$\rho_n=\frac{M_n^{\frac{p-1}{2}}}{\sqrt{A(t_n,u_n)}}.$$
Then $v_n$ satisfies
\begin{equation}\label{eq1568}
\left\{
\begin{array}{ll}
-\Delta v_n=v_n^p+t_n\frac{h(x,M_nv_n,M_n\nabla v_n)}{M_n^p},&y\in\Omega_n,\\
v_n>0,&y\in\Omega_n,\\
v_n=0,&y\in\partial\Omega_n.
\end{array}
\right.
\end{equation}

We claim that, up to a subsequence, we have
$$\lim\limits_{n\rightarrow+\infty}\rho_n=+\infty.$$
To this end, we set $\eta_n=(t_nb)^{\frac{1}{2\alpha}}||\nabla u_n||$. Then $A(t_n,u_n)$ can be rewritten as
$$A(t_n,u_n)=a+\eta_n^{2\alpha}.$$

If $\eta_n$ is bounded by a constant $C$, then, as $n\rightarrow+\infty$, we have
$$\rho_n=\frac{M_n^{\frac{p-1}{2}}}{\sqrt{a+\eta_n^{2\alpha}}}\geq\frac{M_n^{\frac{p-1}{2}}}{\sqrt{a+C^{2\alpha}}}\rightarrow+\infty.$$

If $\eta_n$ is unbounded, then, up to a subsequence, we may assume that
$$\lim\limits_{n\rightarrow+\infty}\eta_n=+\infty.$$
Consequently, we have
$$\sqrt{A(t_n,u_n)}\leq\sqrt{2}\eta_n^{\alpha}.$$
for sufficiently large $n$. By the differential equation in (\ref{eq1}) and $(H_2)$, we can deduce
$$A(t_n,u_n)||\nabla u_n||_2^2=||u_n||_{p+1}^{p+1}+u_nh(x,u_n,\nabla u_n)\leq CM_n^{p+1},\ \ \ \mbox{for large enough}\ \ n.$$
Noting that $t_n\leq 1$, we have
$$||\nabla u_n||_2\geq b^{-\frac{1}{2\alpha}}\eta_n.$$
Therefore, there exist a positive constant $C$ such that
$$A(t_n,u_n)||\nabla u_n||_2^2\geq C\eta_n^{2(\alpha+1)}.$$
This lead to
$$\eta_n\leq CM_n^{\frac{p+1}{2(\alpha+1)}}.$$
Consequently, we have
$$\sqrt{A(t_n,u_n)}\leq CM_n^{\frac{\alpha(p+1)}{2(\alpha+1)}}.$$
Noting that $p>2\alpha+1$, we have
$$\rho_n=\frac{M_n^{\frac{p-1}{2}}}{\sqrt{A(t_n,u_n)}}\geq CM_n^{\frac{p-2\alpha-1}{2(\alpha+1)}}\rightarrow+\infty,\ \ \ \mbox{as}\ \ \ n\rightarrow+\infty.$$
At this stage, we have the following two cases to be considered.

{\bf Case I.}\ If $\{\mbox{dist}(x_n,\partial\Omega_n)\}_{n=1}^\infty$ is unbounded, then, up to a subsequence, we may assume that $\mbox{dist}(x_n,\partial\Omega_n)\rightarrow\infty$ and $\Omega_n\rightarrow R^N$ as $n\rightarrow\infty$. Therefore, for any fixed compact domain $K\subset R^N$ there exist $n_K$ such that $K\subset\Omega_n$ when $n>n_K$. Since $0<v_n\leq 1$, we can get from the assumption $(H_2)$ that
$$|v_n^p+t_n\frac{h(x,M_nv_n,M_n\nabla v_n)}{M_n^p}|\leq C,$$
for some universal positive constant $C$. By the standard regularity theory of elliptic equations, we know that there exists a positive constant $C$ independent of $n$ such that for some $\tau\in(0,1)$ there holds
$$||v_n||_{C^{1,\tau}(K)}\leq C,\ \ \ \mbox{for}\ \ \ \ n>n_K.$$
Hence, up to a subsequence, we have $v_n$ converges in $C^1(K)$ to some function $v$. Since $K$ is arbitrary, by a diagonal process, we can choose a subsequence $\{v_{n_k}\}$ of $\{v_n\}$ such that $v_{n_k}\rightarrow v$ in $C^1_{loc}(R^N)$ as $k\rightarrow +\infty$. Taking the assumption $(H_2)$ into account, it is easy to see that $v$ satisfies
\begin{equation}
\left\{\begin{array}{ll}\label{26}
-\Delta v=v^p & \mbox{in}\ \ R^N,\\
v(0)=1,\\
v\geq 0 & \mbox{in}\ \ R^N.
\end{array}
\right.
\end{equation}
This contradicts with the conclusion of Lemma 3.3.

{\bf Case II.}\ If $\{\mbox{dist}(x_n,\partial\Omega_n)\}_{n=1}^\infty$ is bounded, then, up to a subsequence, we may assume that $\mbox{dist}(x_n,\partial\Omega_n)\rightarrow c$ for some constant $c$, and $\Omega_n\rightarrow R^N_{-c}$ as $n\rightarrow\infty$. Similar to the Case I, we can find a nonnegative function $v$ such that
\begin{equation}
\left\{\begin{array}{ll}\label{27}
-\Delta v=v^p & \mbox{in}\ \ R^N_{-c},\\
v(0)=1,\\
u=0 & x_N=-c.
\end{array}
\right.
\end{equation}
This contradicts with the conclusion of Lemma 3.4.

Therefore, there exists a universal positive constant $R$ such that for any solution $u$ of (\ref{eq1}), there hold $||u||_{C(\Omega)}\leq R$.

In the case $1<p<2\alpha+1<2^*$, we can get a lower bound uniformly in $t$ for solutions of problem (\ref{eq1}) by a similar argument to that of Theorem 3.8. However, we can not expect to have a upper bound uniformly in $t$ for the solutions of problem (\ref{eq1}) due to Theorem 2.2 of section 2. Instead, we have

\vskip 0.1cm

{\em {\bf Theorem 3.9.}\ Assume that $1<p<2\alpha+1<2^*$, and that $(A)$ and $(H_2)$ hold with $\lambda<a\lambda_1(\Omega)$. Then for any fixed $0<t_0\leq 1$ there exists a  positive constant $R$ which may depending on $t_0$ such that for any solution \mbox{$u\in C^2(\Omega)\cap C(\overline{\Omega})$} of problem (\ref{eq1}) with $t\geq t_0$, there holds $||u||_{C(\Omega)}\leq R$.}
\vskip 0.1cm
{\bf Proof.}\ If $u$ is a solution of problem (\ref{eq1}), then it is easy to see that $u$ satisfies
\begin{equation}\label{eq159}
A(t,u)||\nabla u||_2^2=||u||_{p+1}^{p+1}+t\int_{\Omega}uh(x,u,\nabla u)dx.
\end{equation}
Taking $A(t,u)\geq t_0b||\nabla u||_2^{2\alpha}$ and $(H_2)$ into account, we can conclude that there is a positive constant $C$ such that
$$t_0b||\nabla u||_2^{2(\alpha+1)}\leq C(1+||u||_{p+1}^{p+1}).$$
It follows from the Sobolev inequality that
$$t_0b||u||_{p+1}^{2(\alpha+1)}\leq C(1+||u||_{p+1}^{p+1}).$$
Since $p<2\alpha+1$, we can deduce from the above inequality that
$$||u||_{p+1}\leq C,$$
with $C$ depending on $t_0$, $b$, $\alpha$ and $\Omega$. Since $1<p<2^*$, it follows from a bootstrap argument that there exists a positive constant $C$ which is depending on $t_0$, $b$, $\alpha$ and $\Omega$ such that
$$||u||_{C(\Omega)}\leq C.$$
This completes the proof of Theorem 3.9.

From Theorem 3.7, 3.8 and 3.9, we have the following result on positive solutions of problem (\ref{1}).
\vskip 0.1cm
{\em {\bf Corollary 3.10.}\ Assume that $0<p<1$, $(A)$ and $(H_1)$, or $1<p<2^*$, $(A)$ and $(H_2)$ with $\lambda<a\lambda_1(\Omega)$. Then there exist  positive constants $0<r<R$ which are independent of solutions of problem (\ref{1}) such that for any positive solution $u\in C^2(\Omega)\cap C(\overline{\Omega})$ of problem (\ref{1}) there holds $r\leq||u||_{C(\Omega)}\leq R$.}

\section*{4 Existence results}

\setcounter{section}{4}

\setcounter{equation}{0}
This section devotes to prove some existence results for problem (\ref{1}). To this end, we always assume that $h(x,s\xi)$ satisfies $(A)$. Let $C(\overline{\Omega})$, $C^1(\overline{\Omega})$ and $C^{1,\tau}(\overline{\Omega})$ be Banach space equipped with the standard norm. Set $C_0(\overline{\Omega})=\{u(x)\in C(\overline{\Omega})\ \|\ u=0\ \ x\in\partial\Omega \}$, and $E=C^1(\overline{\Omega})\cap C_0(\overline{\Omega})$. For $t\in [0,1]$, define $L_t:\ E\rightarrow C_0(\overline{\Omega})$ by
$$L_t(u)=\frac{|u|^{p-1}u+th(x,u,\nabla u)}{a+tb||\nabla u||_2^{2\alpha}}.$$
Let $(-\Delta)^{-1}$ denote the inverse operator of $-\Delta$. Define $K_t$ by
$$K_t(u)=(-\Delta)^{-1}\circ L_t(u).$$
It is well known that $K_t$ maps $E$ into $C^{1,\tau}(\overline{\Omega})$ for some $\tau\in(0,\ 1)$, and hence is compact.

Since $0<p<2^*$, it follows from the regularity theory that any fixed point $u$ of $K_t$ is a classical solution of the problem
\begin{equation}
\left\{\begin{array}{ll}\label{eq40}
-(a+tb||\nabla u||_2^{2\alpha})\Delta u=|u|^{p-1}u+th(x,u,\nabla u), & \mbox{in}\ \ \Omega,\\
u=0 & \mbox{on}\ \ \partial\Omega,\\
\end{array}
\right.
\end{equation}
Therefore, finding positive solution is equivalent to finding positive fixed point of $K_t$. To find positive fixed point of $K_t$, we set
$$
\begin{array}{ll}
B(0,r)=\{u(x)\in E\ \|\ ||u||_E<r\ \}\\
\mathcal{C}=\{u(x)\in E\ \|\ u(x)\geq 0,\ x\in\Omega\}\\
B_r=B(0,r)\cap\mathcal{C}\\
D(r,R)=\{u(x)\in\mathcal{C}\ \|\ r<||u||_E<R\ \}
\end{array}
$$
It follows from the maximum principle that $K_t$ maps $\mathcal{C}$ into itself. Let $\mathcal{I_{\mathcal{C}}}$ denote the topological index in
$\mathcal{C}$(see \cite{sDan}). If $p\in(1,\ 2^*)$, we know from \cite{sFLN, sRu, sFY, sDGZ} that there exist $0<r_0<R_0$ such that for any $0<r\leq r_0$ and $R\geq R_0$, there holds $\mathcal{I_{\mathcal{C}}}(K_0,D(r,R))=-1$. If $p\in(0,\ 1)$, it is well known that the following problem has a unique solution
\begin{equation}
\left\{\begin{array}{ll}\label{eq41}
-a\Delta u=u^p & \mbox{in}\ \ \Omega,\\
u>0 & \mbox{in}\ \ \Omega,\\
u=0 & \mbox{on}\ \ \partial\Omega,\\
\end{array}
\right.
\end{equation}
Consequently, There are positive numbers $r_1<R_1$ such that for any $0<r\leq r_1$ and $R\geq R_1$, we have $\mathcal{I_{\mathcal{C}}}(K_0,D(r,R))\neq 0$. Noting that the estimates obtained in Theorem 3.7 and 3.8 are uniform in $t$, we can find $0<r_2\leq\min\{r_0,\ r_1\}$ and
$R_2\geq\max\{R_0,\ R_1\}$ such that for any $t\in [0,\ 1]$, $K_t$ has no fixed point on $\partial B_{r_2}$ and $\partial B_{R_2}$. Thus, for any
$t\in [0,\ 1]$, we have
$$\mathcal{I_{\mathcal{C}}}(K_t,D(r_2,R_2))=\mathcal{I_{\mathcal{C}}}(K_0,D(r_2,R_2))\neq 0.$$
Hence, for any $t\in [0,\ 1]$, $K_t$ has at least one fixed point in $D(r_2,R_2)$. Consequently, for any $t\in [0,\ 1]$, the problem
\begin{equation}
\left\{\begin{array}{ll}\label{eq411}
-(a+tb||\nabla u||_2^{2\alpha})\Delta u=|u|^{p-1}u +th(x,u,\nabla u), & \mbox{in}\ \ \Omega,\\
u>0 & \mbox{in}\ \ \Omega,\\
u=0 & \mbox{on}\ \ \partial\Omega,\\
\end{array}
\right.
\end{equation}
has at least one solution. Therefore, we have.
\vskip 0.1cm
{\em {\bf Theorem 4.1.} Assume that $0<p<1$, $(A)$ and $(H_1)$, or $2\alpha+1<p<2^*$, $(A)$ and $(H_2)$. Then, for $a,\ b>0$, the following problem has at least one solution .
\begin{equation}
\left\{\begin{array}{ll}\label{eq412}
-(a+b||\nabla u||_2^{2\alpha})\Delta u=|u|^{p-1}u+h(x,u,\nabla u) & \mbox{in}\ \ \Omega,\\
u>0 & \mbox{in}\ \ \Omega,\\
u=0 & \mbox{on}\ \ \partial\Omega,\\
\end{array}
\right.
\end{equation}}


\newpage


\begin{thebibliography}{s2}




\bibitem{sA} A. Azzollini, The elliptic Kirchhoff equation in $R^N$ perturbed by a local nonlinearity, Differential Integral Equations 25(5-6) (2012) 543-554

\bibitem{sAC} C. Alves, F. Corr$\hat{e}$a, T. Ma, Positive solutions for a quasilinear elliptic equation of Kirchhoff type, Comput. Math. Appl. 49 (2005) 85-93


\bibitem{sCK} C. Chen, Y. Kuo, T. Wu, The Nehari manifold for a Kirchhoff type problem involving sign-changing weight functions, J.Differential Equations 250 (2011) 1876-1908

\bibitem{sCL} M. Chipot, B. Lovat, Some remarks on nonlocal elliptic and parabolic problems, Nonlinear Anal. 30 (1997) 4619-4627

\bibitem{sCW}  B. T. Cheng, X. Wu, Existence results of positive solutions of Kirchhoff type problems, Nonlinear Anal. 71 (2009) 4883-4892



\bibitem{sDGZ} Q. Y. Dai, Y. G. Gu, J. Y. Zhu, A priori estimates, existence and nonexistence of positive solutions of generalized mean curvature equations, Nonlinear Anal. 74 (2011) 7126-7136

\bibitem{sDan} E. N. Dancer, Fixed point index calculations and applications, in: M. Maetzu, A. Vignoli (Eds), Topological Nonlinear Analysis; Degree, Singularity and Variations, Birkhauser, Boston, 1995, pp. 303-340

\bibitem{sFLN} D. de Figueiredo, P. L. Lions, R. D. Nussbaum, A priori estimates and existence of positive solutions of semilinear elliptic equations, J. Math. Pures Appl. 61 (1982) 41-63

\bibitem{sFY} D. de Figueiredo, J. Yang, A priori bounds for positive solutions of a nonvariational elliptic systerm, Comm. Partial Differential Equations 26 (2001) 2305-2321

\bibitem{sFig} G. M. Figueiredo, Existence of a positive solution for Kirchhoff problem type with critical growth via truncation argument, J. Math. Anal. Appl. 401(2013) 706-713

\bibitem{sFS} G. M. Figueiredo, J. R. Jr., Multiplicity of solutions for a Kirchhoff equation with subcritical or critical growth, Differ. Int. Equ. 25(2012) 853-868

\bibitem{sGS} B. Gidas, J. Spruck, Apriori bounds for positive solution s of nonlinear elliptic equations, Comm. Partial Differential Equations 6 (1981) 883-901

\bibitem{sHZ} X. M. He, W. M. Zou, Existence and concentration behavior of positive solutions for a Kirchhoff equation in $R^3$, J. Differential Equations 252 (2012) 1813-1834

\bibitem{sJW} J. H. Jin, X. Wu, Infinitely many radial solutions for Kirchhoff-type problems in $R^N$, J.Math.Anal.Appl. 369 (2010) 564-574

\bibitem{sLLT} J. Liu, J. F. Liao, C. L. Tang, Positive solutions for Kirchhoff-type equations with critical exponent in $R^N$, J.Math.Anal.Appl. 429 (2015) 1153-1172

\bibitem{sLG} Z. S. Liu, S. J. Guo, On ground states for the Kirchhoff-type problem with a general
critical nonlinearity, J. Math. Anal. Appl. 426 (2015) 267-287

\bibitem{sLLS} Y. H. Li, F. Y. Li, J. P. Shi, Existence of a positive solution to Kirchhoff type problems without compactness conditions, J. Differential Equations 253 (2012) 2285-2294


\bibitem{sLLS1} Z. P. Liang, F. Y. Li, J. P. Shi, Positive solutions to Kirchhoff type equations with nonlinearity having prescribed asymptotic behavior, Ann. I.H.Poincar$\acute{e}$-AN 31 (2014) 155-167

\bibitem{sDN} D. Naimen, Positive solutions of Kirchhoff type elliptic equations involving a critical Sobolev exponent, Nonlinear Differ. Equ. Appl. 21(2014) 885-914

\bibitem{sNW} J.J. Nie, X. Wu, Existence and multiplicity of non-trivial solutions for Schr$\ddot{o}$inger-Kirchhoff-type equations with radial potential, Nonlinear Anal. 75 (2012) 3470-3479

\bibitem{sM} T.F. Ma, Remarks on an elliptic equation of Kirchhoff type, Nonlinear Anal. 63 (2005) e1967-e1977

\bibitem{sMR} T. Ma, J. Rivera, Positive solutions for a nonlinear nonlocal elliptic transmission problem, Appl. Math. Lett. 16 (2003) 243-248

\bibitem{sMZ} A. Mao, Z. Zhang, Sign-changing and multiple solutions of Kirchhoff type problems without the P.S. condition, Nonlinear Anal. 70 (2009) 1275-1287

\bibitem{sPZ} K. Perera, Z. Zhang, Nontrivial solutions of Kirchhoff-type problems via the Yang index, J. Differential Equations 221 (2006) 246-255


\bibitem{sRu}D. Ruiz, A priori estimates and existence of positive solutions for strongly nonlinear problems, J. Differential Equations 199 (2004) 96-114

\bibitem{sST} J. J. Sun, C. L. Tang, Existence and multiplicity of solutions for Kirchhoff type equations, Nonlinear Anal. 74(4) (2011) 1212-1222

\bibitem{sWTXZ} J. Wang, L. X. Tian, J. X. Xu, F. B. Zhang, Multiplicity and concentration of positive solutions for a Kirchhoff type problem with critical growth, J. Differential Equations 253 (2012) 2314-2351

\bibitem{sX} X. Wu, Existence of nontrivial solutions and high energy solutions for Schr$\ddot{o}$inger-Kirchhoff-type equations in $R^N$, Nonlinear Anal. Real World Appl. 12 (2011) 1278-1287

\bibitem{sYZ} Y. Yang, J. Zhang, Nontrivial solutions of a class of nonlocal problems via local linking theory, Appl. Math. Lett. 23 (2010) 377-380

\bibitem{sZP} Z. Zhang, K. Perera, Sign changing solutions of Kirchhoff type problems via invariant sets of descent flow, J. Math. Anal. Appl. 317 (2006) 456-463



\end{thebibliography}
\end{document}